# *Chapter 7*

# Multi-Criteria Decision-Making: Reference-Type Methods

**Zhiyuan Wang[a], Gade Pandu Rangaiah[b, c]**

[a] School of Business, Singapore University of Social Sciences, Singapore 599494, Singapore

[b] Department of Chemical and Biomolecular Engineering, National University of Singapore, Singapore 117585, Singapore

[c] School of Chemical Engineering, Vellore Institute of Technology, Vellore 632014, India

**Abstract**

This chapter describes selected reference-type multi-criteria decision-making (MCDM) methods that rank alternatives by comparing them with one or more reference solutions derived from an alternatives-criteria matrix (ACM). After explaining the idea of constructing positive ideal, negative ideal and/or average reference solutions, the chapter details the algorithmic steps of each method, illustrating them with a common ACM example. The 9 methods covered are: Technique for Order of Preference by Similarity to Ideal Solution (TOPSIS), Gray/Grey Relational Analysis (GRA), VlseKriterijumska Optimizacija I Kompromisno Resenje (VIKOR), Evaluation Based on Distance from Average Solution (EDAS), Multi-attributive Border Approximation Area Comparison (MABAC), Combinative Distance-based Assessment (CODAS), Proximity Indexed Value (PIV), Measurement of Alternatives and Ranking According to Compromise Solution (MARCOS) and Preference Ranking on the Basis of Ideal-average Distance (PROBID). The advantages (e.g., computational simplicity) and limitations (e.g., susceptibility to rank reversal) of each method are discussed. A consolidated summary highlights how the different treatments of reference solutions can ultimately drive variations in the ranking of alternatives, underscoring the value of applying several methods in practice. By studying this chapter, readers can (1) describe the principles and steps of each reference-type method, (2) implement them on an ACM, and (3) choose an appropriate reference-type method for their decision-making problems.

## 7.1    Overview

Multi-criteria decision-making (MCDM) serves as an effective approach for recommending one or, in some cases, several alternatives from a set of available alternatives, based on values of multiple criteria. As the diverse landscape of MCDM unfolds in this book, this chapter focuses on the reference-type MCDM methods. These methods utilize one or more reference points derived from the alternatives-criteria matrix (ACM) of the given MCDM application, as a pivotal step. Subsequently, the performance scores of each alternative are calculated based on these reference points. Finally, the alternatives are ranked based on their calculated performance scores.

This chapter delves into 9 reference-type MCDM methods, detailing their principles and algorithms. These methods are selected based on either their widespread popularity for various applications or their recent development, allowing the book to remain up to date with the latest advancements in MCDM algorithms. Each method is elucidated through a simple numerical example, ensuring that readers gain a clear understanding of its inner workings.

The following sections in this chapter cover 9 reference-type MCDM methods, in chronological order. Readers need not study them sequentially and can read any of them (independent of others).

Our Microsoft Excel based program, EMCDM, described in Chapter 10, includes TOPSIS, GRA, VIKOR, MABAC, CODAS, and PROBID methods.

The learning outcomes of this chapter are: (1) Describe the principles and algorithms of reference-type MCDM methods covered; (2) Apply any of these methods to ACM of an application; and (3) Select one or more reference-type MCDM methods for solving MCDM problems.

Unless otherwise stated, the notation used throughout this chapter are as follows:
- $i \in \{1,2,\ldots,m\}$ and $j \in \{1,2,\ldots,n\}$, where $m$ is the number of rows (i.e., alternatives) and $n$ is the number of columns (i.e., criteria) in the ACM.
- $f_{ij}$ is the value of the $j^{th}$ criterion for the $i^{th}$ alternative in the original ACM.

- $F_{ij}$ is the normalized value of $f_{ij}$.
- $w_j$ is the weight assigned to the $j^{th}$ criterion, and $\sum_{j=1}^{n} w_j = 1$.
- $v_{ij}$ is the value of the $j^{th}$ criterion for the $i^{th}$ alternative in the weighted normalized ACM.

Additionally, the ACM presented in Table 7.1 is employed to walk through the steps of the reference-type MCDM methods covered in this chapter. This ACM is made up of 5 alternatives (i.e., $m = 5$ and $i \in \{1,2,3,4,5\}$) and 3 criteria (i.e., $n = 3$ and $j \in \{1,2,3\}$). Out of the 3 criteria (C1, C2, and C3), C1 and C3 are benefit criteria to be maximized, while C2 is a cost criterion to be minimized. For illustration, the 3 criteria are assigned weights as: $w_1 = 0.25$, $w_2 = 0.4$, and $w_3 = 0.35$. All the 5 alternatives (A1, A2, A3, A4, and A5) are non-dominated, which means that the improvement in one criterion will inevitably result in degradation in at least one other criterion (i.e., trade-off among the criteria),

Table 7.1: ACM employed for illustration of reference-type MCDM methods (C1 and C3 for maximization; C2 for minimization).

| Alternatives | C1 (Max) | C2 (Min) | C3 (Max) |
|---|---|---|---|
| A1 | 0.185 | 2.33 | 454 |
| A2 | 0.317 | 1.08 | 298 |
| A3 | 0.555 | 6.45 | 174 |
| A4 | 0.731 | 8.88 | 849 |
| A5 | 0.948 | 7.39 | 517 |
| Weight | 0.25 | 0.4 | 0.35 |

**7.2　Technique for Order of Preference by Similarity to Ideal Solution (TOPSIS)**

The TOPSIS method, introduced by Hwang and Yoon (1981), is centered on the core principle that the top-ranked alternative is the one farthest from the negative ideal solution (NIS) and closest to the positive ideal solution (PIS). In this context, the PIS is composed of the best value for each criterion, that is, the largest value for the criteria to be maximized and the smallest value for the criteria to be minimized (Nabavi et al., 2024). On the contrary, the NIS consists of the worst value for each criterion: the smallest value for maximization criteria and the largest value for minimization criteria. The procedure of the TOPSIS method consists of a series of steps given below; each step is accompanied by numerical calculations using the ACM in Table 7.1 for clear illustration.

**Step 1.** Normalize the original ACM with $m$ rows (i.e., alternatives) and $n$ columns (i.e., criteria) by using Vector normalization method (as in Hwang and Yoon, 1981).

$$F_{ij} = \frac{f_{ij}}{\sqrt{\sum_{k=1}^{m} f_{kj}^2}} \tag{7.1}$$

Here, $f_{kj}$ is the value of the $j^{th}$ criterion for the $k^{th}$ alternative in the original ACM.

*Numerical Calculations (using ACM from Table 7.1):*

For instance, for $i = 2$ and $j = 1$:

$$F_{21} = \frac{f_{21}}{\sqrt{\sum_{k=1}^{5} f_{k1}^2}} = \frac{0.317}{\sqrt{0.185^2 + 0.317^2 + 0.555^2 + 0.731^2 + 0.948^2}} = \frac{0.317}{1.3696} = 0.2315$$

Similar calculations are performed for all other values in the ACM across the 5 alternatives and 3 criteria. The complete results, forming the normalized ACM, are shown in Table 7.2.

Table 7.2: Normalized ACM for TOPSIS walkthrough.

| Alternatives | C1 | C2 | C3 |
|---|---|---|---|
| A1 | 0.1351 | 0.1729 | 0.3962 |
| A2 | 0.2315 | 0.0801 | 0.26 |
| A3 | 0.4052 | 0.4785 | 0.1518 |
| A4 | 0.5337 | 0.6588 | 0.7408 |
| A5 | 0.6922 | 0.5483 | 0.4511 |

**Step 2.** Construct the weighted normalized ACM by multiplying the values of each criterion with its assigned weight, $w_j$.

$$v_{ij} = F_{ij} \times w_j \tag{7.2}$$

*Numerical Calculations:*

As aforementioned, the assigned weights for the 3 criteria are: $w_1 = 0.25$, $w_2 = 0.4$, and $w_3 = 0.35$. Then, for $i = 2$ and $j = 1$:

$$v_{21} = F_{21} \times w_1 = 0.2315 \times 0.25 = 0.0579$$

Likewise, calculations are carried out for all other values in the normalized ACM, using their respective assigned weights. The results, which constitute the weighted normalized ACM, are displayed in Table 7.3.

Table 7.3: Weighted normalized ACM for TOPSIS walkthrough.

| Alternatives | C1 | C2 | C3 |
|---|---|---|---|
| A1 | 0.0338 | 0.0691 | 0.1387 |
| A2 | 0.0579 | 0.0321 | 0.0910 |
| A3 | 0.1013 | 0.1914 | 0.0531 |
| A4 | 0.1334 | 0.2635 | 0.2593 |
| A5 | 0.1730 | 0.2193 | 0.1579 |

**Step 3.** From the weighted normalized ACM, determine the PIS, $A^+$ and NIS, $A^-$ by the following definitions.

$$A^+ = \left\{ \left(\max_i(v_{ij}) \,\middle|\, j \in J_{max}\right), \left(\min_i(v_{ij}) \,\middle|\, j \in J_{min}\right) \right\} \quad (7.3)$$
$$= \{v_1^+, v_2^+, \ldots, v_j^+, \ldots, v_n^+\}$$

$$A^- = \left\{ \left(\min_i(v_{ij}) \,\middle|\, j \in J_{max}\right), \left(\max_i(v_{ij}) \,\middle|\, j \in J_{min}\right) \right\} \quad (7.4)$$
$$= \{v_1^-, v_2^-, \ldots, v_j^-, \ldots, v_n^-\}$$

In the above equations, $J_{max}$ is the set of maximization criteria and $J_{min}$ is the set of minimization criteria, from $\{1, 2, \ldots, n\}$.

*Numerical Calculations:*

Given that criteria C1 and C3 are to be maximized, we have $J_{max} = \{C1, C3\}$; and, as criterion C2 requires minimization, we have $J_{min} = \{C2\}$.

For instance, for $j = 1$:

$$v_1^+ = \max_i(v_{i1}) = 0.1730$$
$$v_1^- = \min_i(v_{i1}) = 0.0338$$

For $j = 2$:

$$v_2^+ = \min_i(v_{i2}) = 0.0321$$
$$v_2^- = \max_i(v_{i2}) = 0.2635$$

The PIS is found to be: $A^+ = \{v_1^+, v_2^+, v_3^+\} = \{0.1730, 0.0321, 0.2593\}$. The NIS is found to be: $A^- = \{v_1^-, v_2^-, v_3^-\} = \{0.0338, 0.2635, 0.0531\}$. Both these can be easily verified by inspecting values in Table 7.3.

**Step 4.** Compute the Euclidean distances of each alternative to both PIS and NIS, respectively, by:

$$\text{Distance to PIS}, S_{i+} = \sqrt{\sum_{j=1}^{n} (v_{ij} - v_j^+)^2} \quad (7.5)$$

$$\text{Distance to NIS}, S_{i-} = \sqrt{\sum_{j=1}^{n}(v_{ij} - v_j^-)^2} \tag{7.6}$$

*Numerical Calculations*:

For instance, for $i = 1$:

$$S_{1+} = \sqrt{\sum_{j=1}^{3}(v_{1j} - v_j^+)^2}$$

$$= \sqrt{(0.0338 - 0.1730)^2 + (0.0691 - 0.0321)^2 + (0.1387 - 0.2593)^2} = 0.1880$$

$$S_{1-} = \sqrt{\sum_{j=1}^{3}(v_{1j} - v_j^-)^2}$$

$$= \sqrt{(0.0338 - 0.0338)^2 + (0.0691 - 0.2635)^2 + (0.1387 - 0.0531)^2} = 0.2124$$

After calculating the Euclidean distances of each alternative to both PIS and NIS, the results are summarized in the first 3 columns of Table 7.4.

Table 7.4: Euclidean distances of each alternative to both PIS and NIS, as well as performance scores ($P_i$), for TOPSIS walkthrough.

| Alternatives | $S_{i+}$ | $S_{i-}$ | $P_i$ |
|---|---|---|---|
| A1 | 0.1880 | 0.2124 | 0.5305 |
| A2 | 0.2039 | 0.2358 | 0.5362 |
| A3 | 0.2703 | 0.0988 | 0.2677 |
| A4 | 0.2348 | 0.2290 | 0.4937 |
| A5 | 0.2130 | 0.1798 | 0.4578 |

**Step 5.** Compute the performance score ($P_i$) of each alternative as follows:

$$P_i = \frac{S_{i-}}{S_{i-} + S_{i+}} \tag{7.7}$$

Finally, the alternative with the largest $P_i$ is top-ranked and recommended to decision-maker.

*Numerical Calculations*:

For instance, for $i = 1$:

$$P_1 = \frac{S_{1-}}{S_{1-} + S_{1+}} = \frac{0.2124}{0.2124 + 0.1880} = 0.5305$$

The calculated performance scores (presented in the last column of Table 7.4) for all alternatives indicate that A2 > A1 > A4 > A5 > A3, with A2 having the largest score of 0.5362, making it the top-ranked alternative by TOPSIS.

The **advantages of TOPSIS** are as follows. (1) As one of the earliest MCDM methods, it is widely adopted in many academic and practical applications, including supply chain management, design, engineering, business, health and safety (Behzadian et al., 2012). (2) Using both PIS and NIS, it is intuitive for decision-makers to accept that the top-ranked solution should be closer to positive ideal while farther away from negative ideal. (3) Steps in the method are relatively straightforward to understand and implement.

The **limitations of TOPSIS** are as follows. (1) It uses the square-based Euclidean distance, which may overly penalize large deviations in one criterion. (2) It is more sensitive to outliers, i.e., one single extreme alternative can change the PIS and NIS. (3) The addition or removal of alternatives can possibly change the PIS and NIS, leading to rank reversal phenomenon (that is, relative rankings of existing alternatives change due to the updated reference points).

## 7.3    Gray/Grey Relational Analysis (GRA)

The GRA method is founded on the gray system theory proposed by Deng (1982). This theory serves as a powerful analytical tool, particularly in situations where information is incomplete or uncertain (Scarlat, 2015). The core principle of GRA is to evaluate alternatives by measuring their proximity to an ideal reference solution through gray relational coefficients (GRCs). The alternative that is closer to the reference solution receives higher ranking (Wang et al., 2023).

There are multiple variants of GRA in the MCDM literature. The specific variant referenced here, as outlined in Song and Jamalipour (2005) and Martinez-Morales et al. (2010), is notable for its autonomy from the predefined weights for criteria. A fundamental step in this GRA approach involves identifying the PIS, which serves as the reference point. The method then computes the distances between each alternative and this PIS to derive the GRC, which measures the similarity of each alternative to the reference point. This coefficient is then used to rank the alternatives, with higher values indicating greater similarity to the ideal solution. The GRA method follows a structured sequence of steps.

**Step 1.** Normalize the original ACM with $m$ rows (i.e., alternatives) and $n$ columns (i.e., criteria) by using the Max-Min normalization method.

$$F_{ij} = \frac{f_{ij} - \min_{k \in \{1,2,\ldots,m\}} f_{kj}}{\max_{k \in \{1,2,\ldots,m\}} f_{kj} - \min_{k \in \{1,2,\ldots,m\}} f_{kj}} \quad \text{for maximization criterion} \quad (7.8)$$

$$F_{ij} = \frac{\max_{k \in \{1,2,\ldots,m\}} f_{kj} - f_{ij}}{\max_{k \in \{1,2,\ldots,m\}} f_{kj} - \min_{k \in \{1,2,\ldots,m\}} f_{kj}} \quad \text{for minimization criterion} \quad (7.9)$$

Note that the normalization equation employed for a minimization criterion reverses its optimization direction, changing from minimizing $f_{ij}$ to maximizing $F_{ij}$ due to the introduction of a negative sign to $f_{ij}$ in the numerator. In addition, all values of $F_{ij}$ are in the range of [0, 1].

*Numerical Calculations (using ACM from Table 7.1):*

For instance, for $i = 2$ and $j = 1$:

$$F_{21} = \frac{0.317 - 0.185}{0.948 - 0.185} = 0.1730$$

Similar calculations are performed for all other values in the ACM, and the complete normalized ACM results are presented in Table 7.5.

Table 7.5: Normalized ACM for GRA walkthrough.

| Alternatives | C1 | C2 | C3 |
|---|---|---|---|
| A1 | 0.0000 | 0.8397 | 0.4148 |
| A2 | 0.1730 | 1.0000 | 0.1837 |
| A3 | 0.4849 | 0.3115 | 0.0000 |
| A4 | 0.7156 | 0.0000 | 1.0000 |
| A5 | 1.0000 | 0.1910 | 0.5081 |

**Step 2.** Find the reference point (i.e., PIS, $A^+$) from the normalized ACM.

$$\begin{aligned} A^+ &= \{\max_i(F_{ij})\} \\ &= \{F_1^+, F_2^+, \ldots, F_j^+, \ldots, F_n^+\} \end{aligned} \quad (7.10)$$

This equation to determine $A^+$ differs from that in TOPSIS due to the application of the Max-Min normalization method in the previous step, which converts minimization criteria into maximization type.

*Numerical Calculations:*

For instance, for $j = 1$:

$$F_1^+ = \max_i(F_{i1}) = 1.0000$$

The reference point is found to be $A^+ = \{F_1^+, F_2^+, F_3^+\} = \{1.0000, 1.0000, 1.0000\}$. This can be readily verified by inspecting values in Table 7.5.

**Step 3.** Construct the difference matrix (Δ) between $F_j^+$ and $F_{ij}$.

$$\Delta_{ij} = |F_j^+ - F_{ij}| \qquad (7.11)$$

*Numerical Calculations:*

For instance, for $i = 2$ and $j = 1$:

$$\Delta_{21} = |F_1^+ - F_{21}| = |1.0000 - 0.1730| = 0.8270$$

Similar calculations are performed for all other values in the normalized ACM, and the complete difference matrix is shown in Table 7.6.

Table 7.6: Difference matrix (Δ) for GRA walkthrough.

| Alternatives | C1 | C2 | C3 |
|---|---|---|---|
| A1 | 1.0000 | 0.1603 | 0.5852 |
| A2 | 0.8270 | 0.0000 | 0.8163 |
| A3 | 0.5151 | 0.6885 | 1.0000 |
| A4 | 0.2844 | 1.0000 | 0.0000 |
| A5 | 0.0000 | 0.8090 | 0.4919 |

**Step 3.** Construct the GRC matrix.

$$GRC_{ij} = \frac{\Delta_{min} + \Delta_{max}}{\Delta_{ij} + \Delta_{max}} \qquad (7.12)$$

Where $\Delta_{max} = \max_{i,j}(\Delta_{ij})$ and $\Delta_{min} = \min_{i,j}(\Delta_{ij})$, which are the global maximum and minimum from the entire difference matrix (Δ), respectively. Note that, by using the Max-Min normalization method, the above GRC values in any application will be between 0.5 and 1.0 (e.g., see Table 7.7).

*Numerical Calculations:*

$$\Delta_{max} = \max_{i,j}(\Delta_{ij}) = 1$$

$$\Delta_{min} = \min_{i,j}(\Delta_{ij}) = 0$$

For instance, for $i = 2$ and $j = 1$:

$$GRC_{21} = \frac{\Delta_{min} + \Delta_{max}}{\Delta_{21} + \Delta_{max}} = \frac{0 + 1}{0.8270 + 1} = 0.5473$$

The complete GRC matrix is shown in Table 7.7.

Table 7.7: GRC matrix for GRA walkthrough.

| Alternatives | C1 | C2 | C3 |
|---|---|---|---|
| A1 | 0.5000 | 0.8619 | 0.6308 |
| A2 | 0.5473 | 1.0000 | 0.5506 |
| A3 | 0.6600 | 0.5923 | 0.5000 |
| A4 | 0.7786 | 0.5000 | 1.0000 |
| A5 | 1.0000 | 0.5528 | 0.6703 |

**Step 4.** Compute the performance score ($P_i$) of each alternative, and the alternative with the largest $P_i$ is top-ranked and recommended to decision-maker.

$$P_i = \frac{1}{n}\sum_{j=1}^{n} GRC_{ij} \tag{7.13}$$

*Numerical Calculations:*

For instance, for $i = 2$:

$$P_2 = \frac{1}{3}\sum_{j=1}^{3} GRC_{2j} = \frac{1}{3}(0.5473 + 1.0000 + 0.5506) = 0.6993$$

Similarly, the performance scores for all alternatives are calculated as $P_1 = 0.6642$, $P_2 = 0.6993$, $P_3 = 0.5841$, $P_4 = 0.7595$, $P_5 = 0.7410$. Accordingly, the ranking is A4 ≻ A5 ≻ A2 ≻ A1 ≻ A3, with A4 being the top-ranked alternative by GRA.

In another widely used GRA variant (Zhang et al., 2011; Kundakçı, 2016), Equation 7.12 in Step 3 is modified as:

$$GRC_{ij} = \frac{\Delta_{min} + \zeta\Delta_{max}}{\Delta_{ij} + \zeta\Delta_{max}} \tag{7.14}$$

Here, $\zeta \in [0,1]$ is the distinguishing coefficient, typically set to 0.5 in the literature, providing moderate distinguishing effects and good stability (Özçelik and Öztürk, 2014). Additionally, if the criteria are assigned different weights, Equation 7.13 in Step 4 can be adjusted accordingly, as follows:

$$P_i = \sum_{j=1}^{n} w_j GRC_{ij} \tag{7.15}$$

The **advantages of GRA** are as follows. (1) The steps in GRA are conceptually simple and intuitive, where bigger GRC value means an alternative is more like the positive ideal reference. (2) Some GRA variants, such as the one outlined in Martinez-Morales et al. (2010), can function without requiring explicit criteria weights from decision-makers, reducing the subjectivity. (3) GRA originates from gray system theory (Deng, 1982), which was designed to handle limited, incomplete, or uncertain data. In principle, GRA can be readily adapted to missing or imprecise values. (4) By introducing the distinguishing coefficient $\zeta \in [0,1]$, sensitivity of GRA to differences ($\Delta$'s) can be controlled.

The **limitations of GRA** are as follows. (1) Similar to other reference-type methods, GRA is also susceptible to rank reversal, where the ranking of alternatives may change when an alternative is added or removed. (2) GRA uses only the PIS as the reference point, whereas some decision-making problems might benefit from considering the distances to both the best and the worst (i.e., PIS and NIS). (3) Although the distinguishing coefficient $\zeta$ provides flexibility in a GRA variant, deciding its value (e.g., 0.3, 0.5, 0.7) is arbitrary, and the final ranking is sensitive to that choice. There is no fully objective way to set $\zeta$ value, and so it becomes another subjective input from the decision-makers.

## 7.4 VlseKriterijumska Optimizacija I Kompromisno Resenje (VIKOR)

The VIKOR (Serbian phrase standing for multi-criteria optimization and compromise solution) method was introduced by Opricovic (1998). It ranks alternatives based on their closeness to the PIS while balancing group benefit (majority rule) and individual regret. The method calculates two key measures: the group utility measure (S), which represents the overall closeness of an alternative to the PIS by aggregating its weighted deviations across all criteria; and the individual regret measure (R), which accounts for the worst-performing criterion for each alternative by considering its maximum deviation from the ideal value. An alternative is ranked higher if it is closer to the PIS (i.e., lower S measure) and has a lower R measure. These two measures are then combined into the VIKOR index (Q), which determines the final ranking of the alternatives. The sequential steps in the VIKOR method are as follows.

**Step 1.** For each criterion, determine its best ($f_j^+$) and worst ($f_j^-$) values in the original ACM with $m$ rows (i.e., alternatives) and $n$ columns (i.e., criteria).

If the $j^{th}$ criterion is for maximization: $f_j^+ = \max_i f_{ij}$ and $f_j^- = \min_i f_{ij}$ \hfill (7.16a)

If the $j^{th}$ criterion is for minimization: $f_j^+ = \min_i f_{ij}$ and $f_j^- = \max_i f_{ij}$ \hfill (7.16b)

*Numerical Calculations (using ACM from Table 7.1):*

For instance, for $j = 1$, since it is a maximization criterion:
$$f_1^+ = \max_i f_{i1} = 0.948$$
$$f_1^- = \min_i f_{i1} = 0.185$$

Likewise, the best and worst values for the remaining criteria are determined as: $f_2^+ = 1.08$, $f_2^- = 8.88$, $f_3^+ = 849$, $f_3^- = 174$.

**Step 2.** Normalize the ACM by using the following equation, where $F_{ij}$ denotes the normalized deviation of the $i^{th}$ alternative to the PIS at the $j^{th}$ criterion.

$$F_{ij} = \frac{f_j^+ - f_{ij}}{f_j^+ - f_j^-} \tag{7.17}$$

This equation is essentially a slight variation of the conventional Max-Min normalization. By substituting the expressions for $f_j^+$ and $f_j^-$, it becomes evident that the optimization direction for maximization criteria is reversed. This transformation changes the maximizing $f_{ij}$ to minimizing $F_{ij}$, due to the negative sign introduced before $f_{ij}$ in the numerator.

*Numerical Calculations:*

For instance, for $i = 1$ and $j = 2$:
$$F_{12} = \frac{f_2^+ - f_{12}}{f_2^+ - f_2^-} = \frac{1.08 - 2.33}{1.08 - 8.88} = 0.1603$$

Likewise, calculations are carried out for all other values in the ACM. The complete results, composing the weighted ACM, are shown in Table 7.8. Notice that the range of each (normalized) criterion is from 0.0 to 1. Further, each normalized value in Table 7.8 = (1.0 – corresponding normalized value in Table 7.5).

Table 7.8: Normalized ACM for VIKOR walkthrough.

| Alternatives | C1 | C2 | C3 |
|---|---|---|---|
| A1 | 1.0000 | 0.1603 | 0.5852 |
| A2 | 0.8270 | 0.0000 | 0.8163 |
| A3 | 0.5151 | 0.6885 | 1.0000 |
| A4 | 0.2844 | 1.0000 | 0.0000 |
| A5 | 0.0000 | 0.8090 | 0.4919 |

**Step 3.** Compute the group utility measure ($S_i$), which is the weighted sum of normalized deviations of the $i^{th}$ alternative to the PIS (in Table 7.8) across all criteria. Further, compute the individual regret measure ($R_i$), which is the maximum normalized deviation of the $i^{th}$ alternative to the PIS among all criteria, ensuring that the method also considers the most unfavorable criterion for each alternative.

$$S_i = \sum_{j=1}^{n} w_j F_{ij} \tag{7.18}$$

$$R_i = \max_j (w_j F_{ij}) \tag{7.19}$$

*Numerical Calculations:*

For instance, for $i = 1$:

$$S_1 = \sum_{j=1}^{n} w_j F_{1j} = 0.25 \times 1.0000 + 0.4 \times 0.1603 + 0.35 \times 0.5852 = 0.5189$$

$$R_1 = \max_j (w_j F_{1j}) = max(0.25 \times 1.0000, 0.4 \times 0.1603, 0.35 \times 0.5852) = 0.2500$$

Similarly, the $S_i$ and $R_i$ values for all alternatives can be calculated; they are presented in **Error! Reference source not found.**.

Table 7.9: Group utility ($S_i$), individual regret ($R_i$), and VIKOR index ($Q_i$) for VIKOR walkthrough.

| Alternatives | $S_i$ | $R_i$ | $Q_i$ |
| --- | --- | --- | --- |
| A1 | 0.5189 | 0.2500 | 0.0845 |
| A2 | 0.4925 | 0.2857 | 0.1567 |
| A3 | 0.7542 | 0.3500 | 0.8333 |
| A4 | 0.4711 | 0.4000 | 0.5000 |
| A5 | 0.4957 | 0.3236 | 0.2888 |

**Step 4.** Compute the VIKOR index ($Q_i$) of each alternative, and the alternative with the lowest $Q$ value is top-ranked and recommended to decision-makers,

$$Q_i = \gamma \left( \frac{S_i - S^+}{S^- - S^+} \right) + (1 - \gamma) \left( \frac{R_i - R^+}{R^- - R^+} \right) \tag{7.20}$$

Here, $S^+ = \min_i S_i$, $S^- = \max_i S_i$, $R^+ = \min_i R_i$, and $R^- = \max_i R_i$; the terms $\frac{S_i - S^+}{S^- - S^+}$ and $\frac{R_i - R^+}{R^- - R^+}$ are essentially the normalized $S_i$ and $R_i$, respectively. This normalization preserves their optimization direction, meaning that smaller values remain preferable. The parameter $\gamma \in [0,1]$ determines the relative importance of group utility versus individual regret in decision-making: $\gamma > 0.5$ favors majority rule, emphasizing overall group benefit (utility); $\gamma = 0.5$ represents a consensus; and $\gamma < 0.5$ gives more weight to individual regret, acting as a veto mechanism. In general, $\gamma$ is set to 0.5 in the literature (Liao and Xu, 2013; Sałabun et al., 2020), representing a balanced compromise (consensus) between group utility and individual regret.

*Numerical Calculations:*

For instance, taking $\gamma = 0.5$, for $i = 1$:

$$Q_1 = 0.5\left(\frac{S_1 - S^+}{S^- - S^+}\right) + 0.5\left(\frac{R_1 - R^+}{R^- - R^+}\right)$$

$$= 0.5\left(\frac{0.5189 - 0.4711}{0.7542 - 0.4711}\right) + 0.5\left(\frac{0.2500 - 0.2500}{0.4000 - 0.2500}\right) = 0.0845$$

Likewise, $Q_i$ values for all alternatives are computed and presented in **Error! Reference source not found.**. The best alternative is the one with the lowest $Q$ value. Accordingly, the ranking list based on $Q$ is A1 > A2 > A5 > A4 > A3, with A1 being the top-ranked alternative by VIKOR.

If decision-makers seek a straightforward ranking and are primarily interested in identifying the top-ranked alternative, they can conclude the process at this stage after obtaining the ranking list based on $Q$, which is the common practice in the literature (Shekhovtsov and Sałabun, 2020; Kizielewicz et al., 2023). However, for completeness, it is important to note that the 3 individual ranking lists based on $S$, $R$, and $Q$ (all follow the "smaller is better" ranking principle) can be considered together to determine either a compromise solution or a set of compromise solutions (Opricovic and Tzeng, 2004). The alternative $a^{(1)}$ is proposed as the compromise solution if it is top-ranked based on $Q$ and concurrently satisfies the following two conditions:

**Condition-1.** Acceptable advantage:

$$Q(a^{(2)}) - Q(a^{(1)}) \geq \frac{1}{m-1} \tag{7.21}$$

Where $a^{(2)}$ is the alternative ranked second in the ranking list based on $Q$. This condition ensures that the top-ranked alternative is significantly better than the second-best alternative.

**Condition-2.** Acceptable stability in decision making, that is, the alternative $a^{(1)}$ must also be the best ranked by $S$ or/and $R$. This condition ensures that the best-ranked alternative is stable under different evaluation perspectives, whether based on group benefit, individual regret, or a combination of both.

If one of the above conditions is not satisfied, then a set of compromise solutions is proposed, which is comprised of:
- Alternatives $a^{(1)}$ and $a^{(2)}$, if only Condition-2 is not satisfied, or
- Alternatives $a^{(1)}, a^{(2)}, ..., a^{(d)}$, if Condition-1 is not satisfied; here, $a^{(d)}$ is the $d^{th}$ ranked alternative in $Q$ and is determined by the following relation for the maximum $d$.

$$Q(a^{(d)}) - Q(a^{(1)}) < \frac{1}{d-1} \tag{7.22}$$

*Numerical Calculations:*

Checking **Condition-1**: inspecting column $Q_i$ from **Error! Reference source not found.**, the second ranked alternative $a^{(2)}$ is A2.

$$Q(a^{(2)}) - Q(a^{(1)}) = 0.1567 - 0.0845 = 0.0723$$

$$\frac{1}{m-1} = \frac{1}{5-1} = 0.25$$

Using Equation 7.21, since 0.0723 is not greater than equal to 0.25, Condition-1 is not satisfied. Therefore, a set of compromise solutions is proposed, following Equation 7.22.

For instance, when $d = 3$, the third ranked alternative $a^{(3)}$ is A5.

$$Q(a^{(3)}) - Q(a^{(1)}) = 0.2888 - 0.0845 = 0.2044$$

$$\frac{1}{d-1} = \frac{1}{3-1} = 0.5$$

Since $0.2044 < 0.5$, A5 is included in the compromise solution set.

When $d = 4$, the fourth ranked alternative $a^{(4)}$ is A4.

$$Q(a^{(4)}) - Q(a^{(1)}) = 0.5000 - 0.0845 = 0.4155$$

$$\frac{1}{d-1} = \frac{1}{4-1} = 0.3333$$

Since 0.4155 is not less than 0.3333, A4 is not included in the compromise solution set.

Finally, the complete set of compromise solutions in this example consists of the alternatives {A1, A2, A5}.

The **advantages of VIKOR** are as follows. (1) VIKOR uniquely balances overall distance from the best (i.e., group utility or benefit, $S_i$) and the worst single-criterion shortfall (i.e., individual regret, $R_i$). This suits problems where decision-makers would like to avoid terrible performance on any one criterion and keep overall performance high. (2) The parameter $\gamma$ can be intuitively interpreted as the relative importance of group utility versus individual regret. Decision-makers can adjust it to reflect their specific attitudes (e.g., risk-averse might focus more on $R_i$ and use a smaller $\gamma$ value). (3) The last step of VIKOR can include conditions to ensure that the top-ranked alternative is a meaningful compromise. If the best alternative is not distinctly better than the second best, a tie or set of compromise solutions is recommended to decision-makers. (4) VIKOR has been extensively applied in various fields (Gul et al., 2016).

The **limitations of VIKOR** are as follows. (1) As with other reference-type methods that rely on PIS and/or NIS, adding or removing an alternative can cause rank reversal. (2) The parameter $\gamma$ is somewhat subjective. A small shift in $\gamma$ can possibly change the final ranking significantly. Without a clear rationale for choosing $\gamma$, decision-makers might question the final recommended results. (3) VIKOR does not always recommend a single top-ranked alternative but instead a set of compromise solutions in some cases (e.g., the ACM employed in this chapter). This requires an extra step from decision-makers to determine the final selection, which introduces more subjectivity. (4) VIKOR's compromise acceptability conditions and the interplay of $S_i$, $R_i$, and $Q_i$ can appear more complicated than other simpler MCDM methods, making it particularly challenging for users without a strong grasp of the underlying mathematics.

## 7.5    Evaluation Based on Distance from Average Solution (EDAS)

The EDAS method, introduced by Keshavarz Ghorabaee et al. (2015), evaluates alternatives by calculating their positive and negative distances from the average solution. Unlike methods such as TOPSIS or VIKOR that rely on PIS and NIS, EDAS is centered on how alternatives deviate from the average values across criteria. Positive distance from the average (PDA) reflects the extent to which an alternative exceeds the average, while negative distance from the average (NDA) measures the extent to which an alternative falls below the average. The assessment of alternatives is based on higher PDA values and lower NDA values, which indicates better performance relative to the average solution. Alternatives are ranked after computing their performance scores, which are derived by combining the normalized PDA and NDA values. The EDAS method can be performed in a structured stepwise manner, as follows.

**Step 1.** Determine the average solution, $\bar{A}$ of the original ACM with $m$ rows (i.e., alternatives) and $n$ columns (i.e., criteria).

$$\bar{f}_j = \frac{\sum_{i=1}^{m} f_{ij}}{m}$$

$$\bar{A} = \{\bar{f}_1, \bar{f}_2, \bar{f}_3, \dots, \bar{f}_j, \dots, \bar{f}_n\}$$

(7.23)

*Numerical Calculations (using ACM from Table 7.1):*

For instance, for $j = 1$:

$$\bar{f}_1 = \frac{\sum_{i=1}^{5} f_{i1}}{5} = \frac{0.185 + 0.317 + 0.555 + 0.731 + 0.948}{5} = 0.5472$$

Likewise, calculations are carried out for the other criteria. The average solution is found to be: $\bar{A} = \{\bar{f}_1, \bar{f}_2, \bar{f}_3\} = \{0.5472, 5.226, 458.4\}$

**Step 2.** Calculate the PDA and NDA matrices based on the determined average solution.

If the $j^{th}$ criterion is for maximization:

$$PDA_{ij} = \frac{max\left(0, (f_{ij} - \bar{f}_j)\right)}{\bar{f}_j}$$

$$NDA_{ij} = \frac{max\left(0, (\bar{f}_j - f_{ij})\right)}{\bar{f}_j}$$

(7.24a)

If the $j^{th}$ criterion is for minimization:

$$PDA_{ij} = \frac{max\left(0, (\bar{f}_j - f_{ij})\right)}{\bar{f}_j}$$

$$NDA_{ij} = \frac{max\left(0, (f_{ij} - \bar{f}_j)\right)}{\bar{f}_j}$$

(7.24b)

Here, $PDA_{ij}$ and $NDA_{ij}$ represent the positive distance (indicating better-than-average performance) and negative distance (indicating worse-than-average performance) of the $i^{th}$ alternative from the average solution with respect to the $j^{th}$ criterion, respectively.

*Numerical Calculations:*

For instance, for $i = 5$ and $j = 1$, since it is a maximization criterion:

$$PDA_{51} = \frac{max\left(0, (f_{51} - \bar{f}_1)\right)}{\bar{f}_1} = \frac{max(0, (0.948 - 0.5472))}{0.5472} = 0.7325$$

$$NDA_{51} = \frac{max\left(0, (\bar{f}_1 - f_{51})\right)}{\bar{f}_1} = \frac{max(0, (0.5472 - 0.948))}{0.5472} = 0$$

For $i = 3$ and $j = 2$, since it is a minimization criterion:

$$PDA_{32} = \frac{max\left(0, (\bar{f}_2 - f_{32})\right)}{\bar{f}_2} = \frac{max(0, (5.226 - 6.45))}{5.226} = 0$$

$$NDA_{32} = \frac{max\left(0, (f_{32} - \bar{f}_2)\right)}{\bar{f}_2} = \frac{max(0, (6.45 - 5.226))}{5.226} = 0.2342$$

Similar calculations are performed for all other values in the ACM, and the complete PDA and NDA matrices are shown in Table 7.10. Note that, for an alternative and a criterion combination, either PDA or NDA is zero.

Table 7.10: PDA and NDA matrices for EDAS walkthrough.

|  | PDA | | | NDA | | |
|---|---|---|---|---|---|---|
| Alternatives | C1 | C2 | C3 | C1 | C2 | C3 |
| A1 | 0 | 0.5542 | 0 | 0.6619 | 0 | 0.0096 |
| A2 | 0 | 0.7933 | 0 | 0.4207 | 0 | 0.3499 |
| A3 | 0.0143 | 0 | 0 | 0 | 0.2342 | 0.6204 |
| A4 | 0.3359 | 0 | 0.8521 | 0 | 0.6992 | 0 |
| A5 | 0.7325 | 0 | 0.1278 | 0 | 0.4141 | 0 |

**Step 3.** Calculate the weighted sum of PDA and NDA (i.e., SP and SN) for each alternative.

$$SP_i = \sum_{j=1}^{n} w_j PDA_{ij} \tag{7.25}$$

$$SN_i = \sum_{j=1}^{n} w_j NDA_{ij} \tag{7.26}$$

*Numerical Calculations:*

For instance, for $i = 1$:

$$SP_1 = \sum_{j=1}^{3} w_j PDA_{1j} = 0.25 \times 0 + 0.4 \times 0.5542 + 0.35 \times 0 = 0.2217$$

$$SN_1 = \sum_{j=1}^{3} w_j NDA_{1j} = 0.25 \times 0.6619 + 0.4 \times 0 + 0.35 \times 0.0096 = 0.1688$$

Analogously, calculations are carried out for the remaining alternatives. The complete $SP_i$ and $SN_i$ results are shown in **Error! Reference source not found.**.

Table 7.11: Weighted sum ($SP_i$ and $SN_i$) and normalized weighted sum ($NSP_i$ and $NSN_i$) of PDA and NDA for EDAS walkthrough.

| Alternatives | $SP_i$ | $SN_i$ | $NSP_i$ | $NSN_i$ |
|---|---|---|---|---|
| A1 | 0.2217 | 0.1688 | 0.5800 | 0.4568 |
| A2 | 0.3173 | 0.2276 | 0.8303 | 0.2676 |
| A3 | 0.0036 | 0.3108 | 0.0093 | 0.0000 |
| A4 | 0.3822 | 0.2797 | 1.0000 | 0.1002 |
| A5 | 0.2279 | 0.1656 | 0.5962 | 0.4671 |

**Step 4.** Normalize the values of $SP_i$ and $SN_i$ for each alternative.

$$NSP_i = \frac{SP_i}{\max_k(SP_k)} \quad k \in \{1,2,\dots,m\} \tag{7.27}$$

$$NSN_i = 1 - \frac{SN_i}{\max_k(SN_k)} \quad k \in \{1,2,\dots,m\} \tag{7.28}$$

When normalizing $SN_i$, the optimization directions of $SN_i$ and $NSN_i$ become opposite due to the introduction of a negative sign in the equation.

*Numerical Calculations:*

For instance, for $i = 1$:

$$NSP_1 = \frac{SP_1}{\max_k(SP_k)} = \frac{0.2217}{0.3822} = 0.5800$$

$$NSN_1 = 1 - \frac{SN_1}{\max_k(SN_k)} = 1 - \frac{0.1688}{0.3108} = 0.4568$$

The complete $NSP_i$ and $NSN_i$ results are presented in the last 2 columns in **Error! Reference source not found.**.

**Step 5.** Compute the performance score ($P_i$) of each alternative as:

$$P_i = \frac{1}{2}(NSP_i + NSN_i) \tag{7.29}$$

Finally, the alternative with the largest $P_i$ is top-ranked and recommended to decision-maker.

*Numerical Calculations:*

For instance, for $i = 1$:

$$P_1 = \frac{1}{2}(NSP_1 + NSN_1) = 0.5800 + 0.4568 = 0.5184$$

Similarly, the performance scores for all alternatives are calculated as $P_1 = 0.5184$, $P_2 = 0.5490$, $P_3 = 0.0047$, $P_4 = 0.5501$, $P_5 = 0.5316$. Accordingly, the ranking is A4 > A2 > A5 > A1 > A3, with A4 being the top-ranked alternative by EDAS.

The **advantages of EDAS** are as follows. (1) It uses the average of each criterion as its reference solution, which is simple and intuitive. Decision-makers can readily interpret "positive distance" (better-than-average) and "negative distance" (worse-than-average). (2) By using the average solution, it is less affected by extreme values compared to methods that rely on ideal solutions. (3) It is computationally less intensive because of its simple calculations based on deviations from the average solution.

The **limitations of EDAS** are as follows. (1) Though less affected, adding or removing alternatives can still alter the average solution, potentially causing rank reversals. (2) Unlike some methods that explicitly compare each alternative to PIS and NIS, EDAS considers "above or below average". However, in some cases, decision-makers may need absolute benchmarks rather than relative-to-mean comparisons. (3) EDAS calculates the performance score as the average of $NSP_i$ and $NSN_i$ (i.e., equal weights by default). This aggregation might not always capture the true preferences of decision-makers.

## 7.6 Multi-attributive Border Approximation Area Comparison (MABAC)

The MABAC method, introduced by Pamučar and Ćirović (2015), is based on the core principle of finding the performance score for each alternative by determining the sum of its distances from a unique reference solution, known as the border approximation area matrix. This matrix is made up of the geometric mean of each weighted normalized criterion, serving as a midway or baseline level of performance for comparison (Wang et al., 2022). The ranking is determined by evaluating how far each alternative deviates from this baseline, with a greater distance indicating a higher-ranked alternative (Nabavi et al., 2023). The implementation of the MABAC method progresses through the following steps.

**Step 1.** Normalize the original ACM with $m$ rows (i.e., alternatives) and $n$ columns (i.e., criteria) by using the Max-Min normalization method.

$$F_{ij} = \frac{f_{ij} - \min_{k \in \{1,2,\ldots,m\}} f_{kj}}{\max_{k \in \{1,2,\ldots,m\}} f_{kj} - \min_{k \in \{1,2,\ldots,m\}} f_{kj}} \quad \text{for maximization criterion} \quad (7.30)$$

$$F_{ij} = \frac{\max_{k \in \{1,2,\ldots,m\}} f_{kj} - f_{ij}}{\max_{k \in \{1,2,\ldots,m\}} f_{kj} - \min_{k \in \{1,2,\ldots,m\}} f_{kj}} \quad \text{for minimization criterion} \quad (7.31)$$

Notably, the above normalization equation for a minimization criterion reverses its optimization direction, changing from minimizing $f_{ij}$ to maximizing $F_{ij}$ because of the introduction of a negative sign to $f_{ij}$ in the numerator. Additionally, all values of $F_{ij}$ are in the range of [0, 1].

*Numerical Calculations (using ACM from Table 7.1):*

For instance, for $i = 2$ and $j = 1$:

$$F_{21} = \frac{0.317 - 0.185}{0.948 - 0.185} = 0.1730$$

Similar calculations are performed for all other values in the ACM, and the complete normalized ACM results are presented in Table 7.12.

Table 7.12: Normalized ACM for MABAC walkthrough.

| Alternatives | C1 | C2 | C3 |
|---|---|---|---|
| A1 | 0.0000 | 0.8397 | 0.4148 |
| A2 | 0.1730 | 1.0000 | 0.1837 |
| A3 | 0.4849 | 0.3115 | 0.0000 |
| A4 | 0.7156 | 0.0000 | 1.0000 |
| A5 | 1.0000 | 0.1910 | 0.5081 |

**Step 2.** Construct the weighted normalized ACM using the assigned weights, $w_j$.

$$v_{ij} = (1 + F_{ij}) \times w_j \tag{7.32}$$

Note that a constant value of 1 is added to all $F_{ij}$ values because the Max-Min normalization method inherently creates a value of 0 in each normalized criterion. This would cause the multiplication to become 0 in the subsequent step of calculating the geometric mean values. To prevent this, a constant value of 1 is introduced, as in Pamučar and Ćirović (2015).

*Numerical Calculations:*

For instance, for $i = 2$ and $j = 1$:

$$v_{21} = (1 + F_{21}) \times w_1 = (1 + 0.1730) \times 0.25 = 0.2933$$

Likewise, calculations are carried out for all other values in the normalized ACM, using their respective assigned weights. The results constituting the weighted normalized ACM are presented in Table 7.13.

Table 7.13: Weighted normalized ACM for MABAC walkthrough.

| Alternatives | C1 | C2 | C3 |
|---|---|---|---|
| A1 | 0.2500 | 0.7359 | 0.4952 |
| A2 | 0.2933 | 0.8000 | 0.4143 |
| A3 | 0.3712 | 0.5246 | 0.3500 |
| A4 | 0.4289 | 0.4000 | 0.7000 |
| A5 | 0.5000 | 0.4764 | 0.5279 |

**Step 3.** Determine the border approximation area matrix B. Essentially, $b_j$ sets the benchmark above which an alternative is comparatively good (for that criterion) and below which it is comparatively poor.

$$b_j = \left( \prod_{i=1}^{m} v_{i1} \right)^{1/m} \tag{7.33}$$

$$B = \{b_1, b_2, \ldots, b_j, \ldots, b_n\}$$

Note that border approximation area matrix B is a row vector with n columns.

*Numerical Calculations:*

For instance, for $j = 1$:

$$b_1 = \left(\prod_{i=1}^{5} v_{i1}\right)^{1/5} = (0.2500 \times 0.2933 \times 0.3712 \times 0.4289 \times 0.5000)^{1/5} = 0.3575$$

Likewise, the border approximation area matrix B is found to be:

$$B = \{b_1, b_2, b_3\} = \{0.3575, 0.5675, 0.4839\}$$

**Step 4.** Compute the performance score ($P_i$) of each alternative.

$$P_i = \sum_{j=1}^{n}(v_{ij} - b_j) \qquad (7.34)$$

The alternative with the largest $P_i$ is top-ranked and recommended to decision-maker.

*Numerical Calculations:*

For instance, for $i = 2$: $P_2 = \sum_{j=1}^{3}(v_{2j} - b_j)$

$$= (0.2933 - 0.3575) + (0.8000 - 0.5675) + (0.4143 - 0.4839) = 0.0987$$

Similarly, the performance scores for all alternatives are calculated as $P_1 = 0.0722$, $P_2 = 0.0987$, $P_3 = -0.1630$, $P_4 = 0.1201$, $P_5 = 0.0954$. Accordingly, the ranking is A4 > A2 > A5 > A1 > A3, with A4 being the top-ranked alternative by MABAC.

The **advantages of MABAC** are as follows. (1) Its novelty lies in utilizing the border approximation area matrix, derived from the geometric mean, as the reference solution. Alternatives are evaluated based on whether they lie above or below this border across the criteria. (2) The computational procedure is straightforward; after normalization and weighting, the position of each alternative relative to the border is determined through a simple (matrix) subtraction. (3) Unlike methods that rely on PIS and/or NIS, MABAC assesses alternatives relative to the geometric mean–based border approximation, making it less susceptible to the extreme values.

The **limitations of MABAC** are as follows. (1) The concept of the "border" may be abstract for decision-makers accustomed to explicit PIS and/or NIS as reference points. In other words, comparing alternatives to an intermediate baseline may feel less intuitive. (2) Although MABAC is less sensitive to extreme values than some methods, outliers can still distort the geometric mean–based border approximation area matrix, affecting all subsequent calculations. (3) Also, it suffers from rank reversal in the case of addition or removal of alternatives.

## 7.7 Combinative Distance-based Assessment (CODAS)

The CODAS method, developed by Keshavarz Ghorabaee et al. (2016), is underpinned by the core principle of finding the performance score for each alternative by calculating both Euclidean and Taxicab distances from the NIS reference point. The alternative that is farther from the NIS is ranked higher (Wang et al., 2024). If there are two alternatives that are incomparable according to Euclidean distance, then Taxicab distance is used as the secondary measure. The steps in the CODAS method are as follows.

**Step 1.** Normalize the original ACM with $m$ rows (i.e., alternatives) and $n$ columns (i.e., criteria) by using the Max normalization method (as in Keshavarz Ghorabaee et al., 2016).

$$F_{ij} = \frac{f_{ij}}{\max_{k \in \{1,2,\ldots,m\}} f_{kj}} \quad \text{for maximization criterion} \quad (7.35)$$

$$F_{ij} = \frac{\min_{k \in \{1,2,\ldots,m\}} f_{kj}}{f_{ij}} \quad \text{for minimization criterion} \quad (7.36)$$

The above normalization equation for a minimization criterion alters the optimization direction of the criterion from minimizing $f_{ij}$ (i.e., the smaller $f_{ij}$, the better) to maximizing $F_{ij}$ (i.e., the larger $F_{ij}$, the better), as $f_{ij}$ is positioned in the denominator of the equation.

---

*Numerical Calculations (using ACM from Table 7.1):*

For instance, for $i = 2$ and $j = 1$:

$$F_{21} = \frac{0.317}{0.948} = 0.3344$$

Similar calculations are performed for all other values in the ACM, and the complete normalized ACM results are presented in Table 7.14.

Table 7.14: Normalized ACM for CODAS walkthrough.

| Alternatives | C1 | C2 | C3 |
|---|---|---|---|
| A1 | 0.1951 | 0.4635 | 0.5347 |
| A2 | 0.3344 | 1.0000 | 0.3510 |
| A3 | 0.5854 | 0.1674 | 0.2049 |
| A4 | 0.7711 | 0.1216 | 1.0000 |
| A5 | 1.0000 | 0.1461 | 0.6090 |

---

**Step 2.** Construct the weighted normalized ACM by multiplying the values of each criterion with its assigned weight, $w_j$.

$$v_{ij} = F_{ij} \times w_j \tag{7.37}$$

*Numerical Calculations:*

For instance, for $i = 2$ and $j = 1$:

$$v_{21} = F_{21} \times w_1 = 0.3344 \times 0.25 = 0.0836$$

Likewise, calculations are carried out for all other values in the normalized ACM, using their respective assigned weights. The resulting weighted normalized ACM is presented in Table 7.15.

Table 7.15: Weighted normalized ACM for CODAS walkthrough.

| Alternatives | C1 | C2 | C3 |
| --- | --- | --- | --- |
| A1 | 0.0488 | 0.1854 | 0.1872 |
| A2 | 0.0836 | 0.4000 | 0.1229 |
| A3 | 0.1464 | 0.0670 | 0.0717 |
| A4 | 0.1928 | 0.0486 | 0.3500 |
| A5 | 0.2500 | 0.0585 | 0.2131 |

**Step 3.** Determine the NIS, $A^-$ from the weighted normalized ACM.

$$\begin{aligned} A^- &= \left\{\min_i(v_{ij})\right\} \\ &= \{v_1^-, v_2^-, \ldots, v_j^-, \ldots, v_n^-\} \end{aligned} \tag{7.38}$$

The above equation to determine $A^-$ differs from that of TOPSIS due to the application of the Max normalization method, which transforms minimization criteria into maximization type.

*Numerical Calculations:*

For instance, for $j = 1$:

$$v_1^- = \min_i(v_{i1}) = 0.0488$$

The NIS is found to be: $A^- = \{v_1^-, v_2^-, v_3^-\} = \{0.0488, 0.0486, 0.0717\}$, which can be verified by reviewing the weighted normalized values in Table 7.15.

**Step 4.** Compute the Euclidean and Taxicab distances of each alternative to the NIS.

$$E_i = \sqrt{\sum_{j=1}^{n}(v_{ij} - v_j^-)^2} \tag{7.39}$$

$$T_i = \sum_{j=1}^{n} |v_{ij} - v_j^-| \tag{7.40}$$

*Numerical Calculations:*

For instance, for $i = 1$:

$$E_1 = \sqrt{\sum_{j=1}^{3} (v_{1j} - v_j^-)^2} = \sqrt{(0.0488 - 0.0488)^2 + (0.1854 - 0.0486)^2 + (0.1872 - 0.0717)^2}$$

$$= 0.1790$$

$$T_1 = \sum_{j=1}^{3} |v_{1j} - v_j^-| = |0.0488 - 0.0488| + |0.1854 - 0.0486| + |0.1872 - 0.0717| = 0.2522$$

Analogously, it can be computed that $E_2 = 0.3568$ and $T_2 = 0.4373$ for $i = 2$. The complete Euclidean and Taxicab distances are presented in Table 7.16.

Table 7.16: Euclidean and Taxicab distances of each alternative to the negative ideal solutions for CODAS walkthrough.

| Alternatives | Euclidean Distance | Taxicab Distance |
|:---:|:---:|:---:|
| A1 | 0.1790 | 0.2522 |
| A2 | 0.3568 | 0.4373 |
| A3 | 0.0993 | 0.1159 |
| A4 | 0.3133 | 0.4223 |
| A5 | 0.2461 | 0.3524 |

**Step 5.** Construct the relative assessment matrix H.

$$h_{ik} = (E_i - E_k) + \psi(E_i - E_k) \times (T_i - T_k) \qquad i, k \in \{1, 2, \ldots, m\} \tag{7.41}$$

Here, $\tau \in [0.01, 0.05]$ and is selected to be 0.02 (as in Keshavarz Ghorabaee et al., 2016). Furthermore, $\psi(E_i - E_k) = 1$ if $|E_i - E_k| \geq \tau$, and $= 0$ if $|E_i - E_k| < \tau$.

*Numerical Calculations:*

For instance, for $i = 2$ and $k = 1$:

$$|E_2 - E_1| = |0.3568 - 0.1790| = 0.1778 \geq 0.02$$
$$\text{and so } \psi(E_2 - E_1) = 1$$

Then,

$$h_{21} = (E_2 - E_1) + \psi(E_2 - E_1) \times (T_2 - T_1)$$
$$= (0.3568 - 0.1790) + 1 \times (0.4373 - 0.2522) = 0.3629$$

Likewise, other values in the H matrix can be computed. The complete H matrix is presented in Table 7.17.

Table 7.17: H matrix for CODAS walkthrough.

| Alternatives | A1 | A2 | A3 | A4 | A5 |
|---|---|---|---|---|---|
| A1 | 0 | -0.3629 | 0.2160 | -0.3044 | -0.1674 |
| A2 | 0.3629 | 0 | 0.5789 | 0.0585 | 0.1955 |
| A3 | -0.2160 | -0.5789 | 0 | -0.5204 | -0.3834 |
| A4 | 0.3044 | -0.0585 | 0.5204 | 0 | 0.1370 |
| A5 | 0.1674 | -0.1955 | 0.3834 | -0.1370 | 0 |

**Step 6.** Compute the performance score ($P_i$) of each alternative as follows:

$$P_i = \sum_{k=1}^{m} h_{ik} \qquad (7.42)$$

The alternative with the largest $P_i$ is top-ranked and recommended to decision-maker.

*Numerical Calculations:*

For instance, for $i = 2$:

$$P_2 = \sum_{k=1}^{5} h_{2k} = 0.3629 + 0 + 0.5789 + 0.0585 + 0.1955 = 1.1957$$

Similarly, performance scores for all alternatives are calculated as $P_1 = -0.6187$, $P_2 = 1.1957$, $P_3 = -1.6986$, $P_4 = 0.9034$, $P_5 = 0.2183$. Accordingly, the ranking is A2 > A4 > A5 > A1 > A3, with A2 being the top-ranked alternative. As can be seen, $P_i$ can be negative.

The **advantages of CODAS** are as follows. (1) The use of Max normalization method effectively transforms minimization criteria into maximization criteria, which slightly simplifies the subsequent calculations in determining the NIS. (2) It uses both Euclidean and Taxicab distances, improving the ability to distinguish alternatives and providing a more nuanced ranking. (3) The introduction of the novel threshold parameter, $\tau$ allows for an integrated use of Euclidean and Taxicab distances, theoretically enabling greater flexibility in addressing complex decision-making scenarios.

The **limitations of CODAS** are as follows. (1) Choosing an appropriate $\tau$ is non-trivial and may require input from the decision-maker; an improper $\tau$ value can cause erratic or inconsistent rankings. (2) CODAS references only the NIS, which might not align with the preferences of decision-makers who favor anchoring their evaluations to the PIS as well. (3) In the original CODAS method, the threshold function $\psi$ is defined as $\psi(E_i - E_k) = 1$ if $|E_i - E_k| \geq \tau$, and $= 0$ if $|E_i - E_k| < \tau$. However, Chen and Gou (2021) pointed out that this approach is inconsistent with realistic decision-making. They stated that, when the Euclidean distance difference between two alternatives is tiny (i.e., less than the threshold), the Taxicab distance should be incorporated to provide additional differentiation. On the contrary, when the Euclidean distance difference between two alternatives is large (i.e., greater than or equal to the threshold), the Euclidean distance alone should suffice to distinguish between the two alternatives, making the use of Taxicab distance unnecessary.

## 7.8 Proximity Indexed Value (PIV)

The PIV method, introduced by Mufazzal and Muzakkir (2018) to mitigate the rank reversal phenomenon, is based on the core principle of determining the performance score of each alternative by comparing it against the PIS. Alternatives that are closer to the PIS receive higher rankings. The weighted proximity index is determined by calculating the difference between the weighted normalized value ($v_{ij}$) and the best value of each criterion. The PIV method is implemented through the following series of steps.

**Step 1.** Normalize the original ACM with $m$ rows (i.e., alternatives) and $n$ columns (i.e., criteria) by using Vector normalization method. This step does not differ from Step 1 of the TOPSIS method, and Equation 7.1 is equally applicable here. For brevity, this equation and numerical calculations (same as Table 7.2) are not repeated here.

**Step 2.** Construct the weighted normalized ACM using the assigned weights, $w_j$. As with Step 1, Equation 7.2 is also applicable here. The equation and numerical calculations (same as Table 7.3) are omitted here.

**Step 3.** Determine the best value for each criterion ($v_j^+$), forming the PIS. This step is the same as Step 3 of the TOPSIS method and utilizes Equation 7.3. The PIS is also found to be: $A^+ = \{v_1^+, v_2^+, v_3^+\} = \{0.1730, 0.0321, 0.2593\}$.

**Step 4.** Calculate the overall proximity value ($D_i$) of each alternative to the PIS by:

$$D_i = \sum_{j=1}^{n} |v_{ij} - v_j^+| \qquad (7.43)$$

The alternative with the smallest $D_i$ is top-ranked and recommended to the decision-maker. In effect, the PIV method finds the alternative closest to the PIS (as measured by $D_i$).

> *Numerical Calculations:*
>
> For instance, for $i = 1$:
>
> $$D_1 = \sum_{j=1}^{3} |v_{1j} - v_j^+| = |0.0338 - 0.1730| + |0.0691 - 0.0321| + |0.1387 - 0.2593| = 0.2970$$
>
> Similarly, the overall proximity values for all alternatives are calculated as $D_1 = 0.2970$, $D_2 = 0.2835$, $D_3 = 0.4373$, $D_4 = 0.2711$, $D_5 = 0.2887$. Accordingly, the ranking is A4 > A2 > A5 > A1 > A3, with A4 being the top-ranked alternative by PIV.

The **advantages of PIV** are as follows. (1) It focuses on the distance from the PIS only, reducing conceptual (and computational) overhead compared to other reference-type methods. (2) The essential usage of L1 norm (i.e., summation of absolute deviations from the PIS) is intuitive for decision-makers. (3) The concept and implementation of PIV are easy to understand.

The **limitations of PIV** are as follows. (1) By not considering the NIS, it fails to capture how far alternatives are from the worst-case scenario, which could be a critical consideration in risk-averse applications. (2) It is more susceptible to the influence of extreme values. If the identified reference solution is an extreme outlier, it can skew all distances significantly. (3) L1 norm treats all deviations linearly with equal emphasis, potentially under-penalizing big shortfall in one criterion. (4) The default use of the Max normalization method may result in a division by zero error.

## 7.9 Measurement of Alternatives and Ranking According to Compromise Solution (MARCOS)

The MARCOS method, introduced by Stević et al. (2020), is based on the principle of measuring alternatives and ranking them according to their compromise with both the PIS and NIS. By establishing the relationship of each alternative to these two reference points, MARCOS determines utility functions and their aggregation, which reflect the proximity of an alternative to the PIS and its distance from the NIS. An alternative is ranked higher if it is closer to PIS and farther from NIS. According to Stević et al. (2020), this method introduces a novel approach to defining and aggregating utility functions, enhancing the decision-making framework. They also performed comprehensive sensitivity analyses, including variations in weight values, decision-making matrices, and comparisons with other MCDM methods, demonstrating its robustness and validity across different scenarios. The MARCOS method operates through the following step-by-step procedure.

**Step 1.** Define the PIS, $A^+$ and the NIS, $A^-$ :

*Preliminary Draft Manuscript*$$A^+ = \left\{\left(\max_i(f_{ij}) \mid j \in J_{max}\right), \left(\min_i(f_{ij}) \mid j \in J_{min}\right)\right\}$$
$$= \{f_1^+, f_2^+, \ldots, f_j^+, \ldots, f_n^+\} \tag{7.44}$$

$$A^- = \left\{\left(\min_i(f_{ij}) \mid j \in J_{max}\right), \left(\max_i(f_{ij}) \mid j \in J_{min}\right)\right\}$$
$$= \{f_1^-, f_2^-, \ldots, f_j^-, \ldots, f_n^-\} \tag{7.45}$$

In both the above equations, $J_{max}$ is the set of maximization criteria and $J_{min}$ is the set of minimization criteria, both from $\{1, 2, \ldots, n\}$. Then, add both the PIS and NIS to the original ACM (with $m$ alternatives and $n$ criteria), forming an extended ACM, as illustrated below.

*Numerical Calculations (using ACM from Table 7.1):*

For instance, for $j = 1$:

$$f_1^+ = \max_i(f_{i1}) = 0.948$$

$$f_1^- = \min_i(f_{i1}) = 0.185$$

For $j = 2$:

$$f_2^+ = \min_i(f_{i2}) = 1.08$$

$$f_2^- = \max_i(f_{i2}) = 8.88$$

The PIS is found to be: $A^+ = \{f_1^+, f_2^+, f_3^+\} = \{0.948, 1.08, 849\}$. The NIS is found to be: $A^- = \{f_1^-, f_2^-, f_3^-\} = \{0.185, 8.88, 174\}$. The extended ACM is shown in Table 7.18, with $A^+$ and $A^-$ appended as the last two rows.

Table 7.18: Extended ACM for MARCOS walkthrough.

| Alternatives | C1 | C2 | C3 |
| --- | --- | --- | --- |
| A1 | 0.185 | 2.33 | 454 |
| A2 | 0.317 | 1.08 | 298 |
| A3 | 0.555 | 6.45 | 174 |
| A4 | 0.731 | 8.88 | 849 |
| A5 | 0.948 | 7.39 | 517 |
| $A^+$ | 0.948 | 1.08 | 849 |
| $A^-$ | 0.185 | 8.88 | 174 |

**Step 2.** Normalize the extended ACM (with $i \in \{1, 2, \ldots, m+2\}$) by using the Max normalization method (as in Stević et al., 2020).

7-29

$$F_{ij} = \frac{f_{ij}}{\max_{k \in \{1,2,\ldots,m\}} f_{kj}} \quad \text{for maximization criterion} \quad (7.46)$$

$$F_{ij} = \frac{\min_{k \in \{1,2,\ldots,m\}} f_{kj}}{f_{ij}} \quad \text{for minimization criterion} \quad (7.47)$$

Use of the above normalization equation for a minimization criterion alters the optimization direction of the criterion from minimizing $f_{ij}$ (i.e., the smaller $f_{ij}$, the better) to maximizing $F_{ij}$ (i.e., the larger $F_{ij}$, the better), as $f_{ij}$ is in the denominator of the equation.

*Numerical Calculations:*

For instance, for $i = 2$ and $j = 1$:

$$F_{21} = \frac{0.317}{0.948} = 0.3344$$

Similar calculations are performed for all other values in the extended ACM, and the complete normalized extended ACM results are presented in Table 7.19.

Table 7.19: Normalized extended ACM for MARCOS walkthrough.

| Alternatives | C1 | C2 | C3 |
|---|---|---|---|
| A1 | 0.1951 | 0.4635 | 0.5347 |
| A2 | 0.3344 | 1.0000 | 0.3510 |
| A3 | 0.5854 | 0.1674 | 0.2049 |
| A4 | 0.7711 | 0.1216 | 1.0000 |
| A5 | 1.0000 | 0.1461 | 0.6090 |
| $A^+$ | 1.0000 | 1.0000 | 1.0000 |
| $A^-$ | 0.1951 | 0.1216 | 0.2049 |

**Step 3.** Construct the weighted normalized extended ACM (with $i \in \{1,2,\ldots,m+2\}$) by multiplying the values of each criterion with its assigned weight, $w_j$.

$$v_{ij} = F_{ij} \times w_j \quad (7.48)$$

As seen from the first 3 steps, MARCOS initially identifies the ideal solutions, followed by normalization and then applying weight. This sequence essentially leads to the same outcome as the conventional approach of normalizing first, then applying weights, and subsequently determining ideal solutions. This can be observed by comparing the following Table 7.20 with Table 7.15 in the CODAS section (as both MARCOS and CODAS employ Max normalization method).

*Numerical Calculations:*

For instance, for $i = 2$ and $j = 1$:

$$v_{21} = F_{21} \times w_1 = 0.3344 \times 0.25 = 0.0836$$

Likewise, calculations are carried out for all other values in the normalized extended ACM, using their respective assigned weights. The results, which constitute the weighted normalized extended ACM, are presented in Table 7.20.

Table 7.20: Weighted normalized extended ACM for MARCOS walkthrough.

| Alternatives | C1 | C2 | C3 |
| --- | --- | --- | --- |
| A1 | 0.0488 | 0.1854 | 0.1872 |
| A2 | 0.0836 | 0.4000 | 0.1229 |
| A3 | 0.1464 | 0.0670 | 0.0717 |
| A4 | 0.1928 | 0.0486 | 0.3500 |
| A5 | 0.2500 | 0.0585 | 0.2131 |
| $A^+$ | 0.2500 | 0.4000 | 0.3500 |
| $A^-$ | 0.0488 | 0.0486 | 0.0717 |

**Step 4.** Calculate the utility degrees ($K_i^+$ and $K_i^-$) and utility function values ($f(K_i^+)$ and $f(K_i^-)$) of alternatives in relation to the PIS and NIS, respectively, by the following equations.

$$S_i = \sum_{j=1}^{n} v_{ij} \qquad \text{for } i \in \{1, 2, \ldots, m\} \qquad (7.49a)$$

$$S_{A^+} = \sum_{j=1}^{n} v_{ij} \qquad \text{for } i = m + 1 \qquad (7.49b)$$

$$S_{A^-} = \sum_{j=1}^{n} v_{ij} \qquad \text{for } i = m + 2 \qquad (7.49c)$$

$$K_i^+ = \frac{S_i}{S_{A^+}} \qquad \text{for } i \in \{1, 2, \ldots, m\} \qquad (7.50)$$

$$K_i^- = \frac{S_i}{S_{A^-}} \qquad \text{for } i \in \{1, 2, \ldots, m\} \qquad (7.51)$$

Here, $S_{A^+}$ and $S_{A^-}$ represent the total sum of the weighted normalized values across all criteria for the PIS and NIS ($A^+$ and $A^-$), respectively.

$$f(K_i^+) = \frac{K_i^-}{K_i^+ + K_i^-} \qquad \text{for } i \in \{1, 2, \ldots, m\} \qquad (7.52)$$

$$f(K_i^-) = \frac{K_i^+}{K_i^+ + K_i^-} \quad \text{for } i \in \{1, 2, \ldots, m\} \tag{7.53}$$

*Numerical Calculations:*

$$S_{A^+} = 0.2500 + 0.4000 + 0.3500 = 1.0000$$

Note that the above quantity will always be unity due to Max normalization and weights chosen to sum to unity.

$$S_{A^-} = 0.0488 + 0.0486 + 0.0717 = 0.1692$$

For instance, for $i = 1$:

$$S_1 = 0.0488 + 0.1854 + 0.1872 = 0.4214$$

$$K_1^+ = \frac{S_1}{S_{A^+}} = \frac{0.4214}{1} = 0.4214$$

$$K_1^- = \frac{S_1}{S_{A^-}} = \frac{0.4214}{0.1692} = 2.4908$$

$$f(K_1^+) = \frac{K_1^-}{K_1^+ + K_1^-} = \frac{2.4908}{0.4214 + 2.4908} = 0.8553$$

$$f(K_1^-) = \frac{K_1^+}{K_1^+ + K_1^-} = \frac{0.4214}{0.4214 + 2.4908} = 0.1447$$

Analogously, utility degrees ($K_i^+$ and $K_i^-$) and utility function values ($f(K_i^+)$ and $f(K_i^-)$) of the other alternatives are calculated and presented in Table 7.21. Notably, $f(K_i^+)$ and $f(K_i^-)$ values are identical for all alternatives because $f(K_i^+)$ equation can be simplified as follows:

$$f(K_i^+) = \frac{K_i^-}{K_i^+ + K_i^-} = \frac{\frac{S_i}{S_{A^-}}}{\frac{S_i}{S_{A^+}} + \frac{S_i}{S_{A^-}}} = \frac{\frac{1}{S_{A^-}}}{\frac{1}{S_{A^+}} + \frac{1}{S_{A^-}}}$$

This simplified expression depends solely on $S_{A^+}$ and $S_{A^-}$, making utility function calculations invariant to the specific alternative.

Table 7.21: Utility degrees and utility function values for MARCOS walkthrough.

| Alternatives | $S_i$ | $K_i^+$ | $K_i^-$ | $f(K_i^+)$ | $f(K_i^-)$ |
|---|---|---|---|---|---|
| A1 | 0.4214 | 0.4214 | 2.4908 | 0.8553 | 0.1447 |
| A2 | 0.6064 | 0.6064 | 3.5849 | 0.8553 | 0.1447 |
| A3 | 0.2851 | 0.2851 | 1.6851 | 0.8553 | 0.1447 |
| A4 | 0.5914 | 0.5914 | 3.4961 | 0.8553 | 0.1447 |
| A5 | 0.5216 | 0.5216 | 3.0833 | 0.8553 | 0.1447 |
| $A^+$ | 1.0000 | | | | |
| $A^-$ | 0.1692 | | | | |

**Step 5.** Compute the performance score ($P_i$) of each alternative as follows:

$$P_i = \frac{K_i^+ + K_i^-}{1 + \frac{1-f(K_i^+)}{f(K_i^+)} + \frac{1-f(K_i^-)}{f(K_i^-)}} \quad \text{for } i \in \{1,2,\ldots,m\} \quad (7.54)$$

The alternative with the largest $P_i$ is top-ranked and recommended to the decision-maker.

*Numerical Calculations:*

For instance, for $i = 1$:

$$P_1 = \frac{K_1^+ + K_1^-}{1 + \frac{1-f(K_1^+)}{f(K_1^+)} + \frac{1-f(K_1^-)}{f(K_1^-)}} = \frac{0.4214 + 2.4908}{1 + \frac{1 - 0.8553}{0.8553} + \frac{1 - 0.1447}{0.1447}} = 0.4113$$

Similarly, performance scores for all alternatives are calculated as $P_1 = 0.4113$, $P_2 = 0.5920$, $P_3 = 0.2783$, $P_4 = 0.5773$, $P_5 = 0.5091$. Accordingly, the ranking is A2 ≻ A4 ≻ A5 ≻ A1 ≻ A3, with A2 being the top-ranked alternative.

The **advantages of MARCOS** are as follows. (1) By explicitly adding the PIS and NIS to the ACM matrix, the process becomes more transparent, allowing decision-makers to see exactly what (raw) reference points are being used as benchmarks. (2) The method has been successfully applied in diverse fields, such as supplier selection, healthcare, logistics, and manufacturing, showcasing its good adaptability (Demir et al., 2024). (3) It delivers a single, intuitively interpretable indicator, i.e., the utility degree, for every alternative. This scalar is computed from the alternative's simultaneous proximity to the PIS and NIS, giving decision-makers an at-a-glance sense of "how good" each option is while still preserving full mathematical rigor.

The **limitations of MARCOS** are as follows. (1) Compared to simpler methods (e.g., EDAS), MARCOS involves more calculations, particularly in the utility function aggregation step, making it more computationally intensive when handling a large ACM. (2) Decision-makers unfamiliar with utility functions may find it challenging to interpret results without a solid understanding of the underlying mathematics. (3) Like other reference-type methods, MARCOS can exhibit rank reversal when the set of alternatives changes.

## 7.10 Preference Ranking on the Basis of Ideal-average Distance (PROBID)

The PROBID method, proposed by Wang et al. (2021), is based on the core principle of determining the performance score of each alternative by comparing it against a range of reference solutions. This range encompasses a hierarchy of ideal solutions, starting from the most PIS and extending through lesser tiers (e.g., 2nd and 3rd most PIS) down to the m[th] PIS (i.e., the most NIS). Additionally, the PROBID method incorporates an average solution to provide a comprehensive baseline for comparison. Hence, according to Wang et al. (2021), this method is less susceptible to minor data changes, reducing the likelihood of rank reversal. By considering multiple layers of ideal solutions and the average solution, PROBID provides a more refined ranking of alternatives (Baydaş et al., 2024). Besides, a simplified variant, sPROBID is also available in the literature (Nabavi et al., 2025). In PROBID, the alternative that is closer to PIS group, farther from NIS group, and deviates more from the average solutions is ranked higher. The procedure of the PROBID method unfolds in a series of steps as follows.

**Step 1.** Normalize the original ACM with $m$ rows (i.e., alternatives) and $n$ columns (i.e., criteria) by using Vector normalization method. This step does not differ from Step 1 of the TOPSIS method, and Equation 7.1 is equally applicable here. For brevity, this equation and numerical calculations (same as Table 7.2) are not repeated here.

**Step 2.** Construct the weighted normalized ACM using the assigned weights, $w_j$. As with Step 1, Equation 7.2 is also applicable here. This equation and numerical calculations (same as Table 7.3) are omitted here.

**Step 3.** Determine the most PIS ($A_{(1)}$), 2nd PIS ($A_{(2)}$), 3rd PIS ($A_{(3)}$), ..., and m[th] PIS ($A_{(m)}$) (i.e., the most NIS) by:

$$A_{(k)} = \left\{ \left( \underset{i}{Large}(v_{ij}, k) \,\big|\, j \in J_{max} \right), \left( \underset{i}{Small}(v_{ij}, k) \,\big|\, j \in J_{min} \right) \right\}$$
$$= \{v_{(k)1}, v_{(k)2}, v_{(k)3}, \dots, v_{(k)j}, \dots, v_{(k)n}\} \quad (7.55)$$

Here, $J_{max}$ is the set of maximization criteria and $J_{min}$ is the set of minimization criteria; $k \in \{1,2,\dots,m\}$; and $\underset{i}{Large}(v_{ij}, k)$ is the $k^{th}$ largest value in the $j^{th}$ criterion (i.e., $v_{(k)j}$). Then, compute the average solution, $\bar{A}$ as follows.

$$\bar{v}_j = \frac{\sum_{i=1}^{m} v_{ij}}{m}$$
$$\bar{A} = \{\bar{v}_1, \bar{v}_2, \bar{v}_3, \dots, \bar{v}_j, \dots, \bar{v}_n\} \quad (7.56)$$

*Numerical Calculations:*

For instance, for $k = 2$, that is, finding the 2nd best value of each criterion:

$$A_{(2)} = \{v_{(2)1}, v_{(2)2}, v_{(2)3}\} = \{0.1334, 0.0691, 0.1579\}$$

Likewise, the other tiers of ideal solutions and the average solution can be determined. The complete results are shown in Table 7.22.

Table 7.22: Tiers of ideal solutions and the average solution for PROBID walkthrough.

| Ideal Solutions | C1 | C2 | C3 |
|---|---|---|---|
| $A_{(1)}$ | 0.1730 | 0.0321 | 0.2593 |
| $A_{(2)}$ | 0.1334 | 0.0691 | 0.1579 |
| $A_{(3)}$ | 0.1013 | 0.1914 | 0.1387 |
| $A_{(4)}$ | 0.0579 | 0.2193 | 0.0910 |
| $A_{(5)}$ | 0.0338 | 0.2635 | 0.0531 |
| Average | 0.0999 | 0.1551 | 0.1400 |

**Step 4.** Iteratively calculate the Euclidean distance of each alternative to each of the *m* ideal solutions as well as to the average solution.

$$\text{Distance to the } k^{th} \text{ PIS}, S_{i(k)} = \sqrt{\sum_{j=1}^{n}(v_{ij} - v_{(k)j})^2}$$

$$\text{Distance to the average solution}, S_{i(avg)} = \sqrt{\sum_{j=1}^{n}(v_{ij} - \bar{v}_j)^2}$$

(7.57)

*Numerical Calculations:*

For instance, for $i = 1$ and $k = 2$:

$$S_{1(2)} = \sqrt{\sum_{j=1}^{3}(v_{1j} - v_{(2)j})^2}$$

$$= \sqrt{(0.0338 - 0.1334)^2 + (0.0691 - 0.0691)^2 + (0.1387 - 0.1579)^2}$$

$$= 0.1015$$

The complete distances of all alternatives are presented in Table 7.23. Note that the $S_{i(1)}$ and $S_{i(5)}$ are the same as the $S_{i+}$ and $S_{i-}$ of TOPSIS (in Table 7.4), respectively.

Table 7.23: Distances of each alternative to the different ideal solutions and average solution for PROBID walkthrough.

| Alternatives | $S_{i(1)}$ | $S_{i(2)}$ | $S_{i(3)}$ | $S_{i(4)}$ | $S_{i(5)}$ | $S_{i(avg)}$ |
|---|---|---|---|---|---|---|
| A1 | 0.1880 | 0.1015 | 0.1397 | 0.1594 | 0.2124 | 0.1084 |
| A2 | 0.2039 | 0.1075 | 0.1719 | 0.1873 | 0.2358 | 0.1389 |
| A3 | 0.2703 | 0.1642 | 0.0855 | 0.0640 | 0.0988 | 0.0942 |
| A4 | 0.2348 | 0.2192 | 0.1442 | 0.1897 | 0.2290 | 0.1647 |
| A5 | 0.2130 | 0.1553 | 0.0793 | 0.1332 | 0.1798 | 0.0990 |

**Step 5.** Determine the overall positive-ideal and negative-ideal distances

$$S_{i(pos-ideal)} = \begin{cases} \sum_{k=1}^{\frac{m+1}{2}} \frac{1}{k} S_{i(k)} & \text{when } m \text{ is an odd number} \\ \sum_{k=1}^{\frac{m}{2}} \frac{1}{k} S_{i(k)} & \text{when } m \text{ is an even number} \end{cases}$$

(7.58)

$$S_{i(neg-ideal)} = \begin{cases} \sum_{k=\frac{m+1}{2}}^{m} \frac{1}{m-k+1} S_{i(k)} & \text{when } m \text{ is an odd number} \\ \sum_{k=\frac{m}{2}+1}^{m} \frac{1}{m-k+1} S_{i(k)} & \text{when } m \text{ is an even number} \end{cases}$$

Note that, when $m$ is an odd number, $S_{i\left(\frac{m+1}{2}\right)}$ is used for calculating both $S_{i(pos-ideal)}$ and $S_{i(ne-ideal)}$.

*Numerical Calculations:*

For instance, for $i = 1$:

$$S_{1(pos-ideal)} = \sum_{k=1}^{3} \frac{1}{k} S_{1(k)} = 1 \times 0.1880 + \frac{1}{2} \times 0.1015 + \frac{1}{3} \times 0.1397 = 0.2853$$

The complete results for the overall positive-ideal and negative-ideal distances of each alternative are presented in Table 7.24.

Table 7.24: Overall positive-ideal and negative-ideal distances for PROBID walkthrough

| Alternatives | $S_{i(pos-ideal)}$ | $S_{i(neg-ideal)}$ |
|---|---|---|
| A1 | 0.2853 | 0.3386 |
| A2 | 0.3150 | 0.3867 |
| A3 | 0.3809 | 0.1593 |
| A4 | 0.3925 | 0.3719 |
| A5 | 0.3171 | 0.2728 |

**Step 6.** Compute the performance score ($P_i$) of each alternative, and the alternative with the largest $P_i$ is top-ranked and recommended to the decision-maker.

$$R_i = \frac{S_{i(pos-ideal)}}{S_{i(neg-ideal)}}$$

$$P_i = \frac{1}{1 + R_i^2} + S_{i(avg)}$$

(7.59)

*Numerical Calculations:*

For instance, for $i = 1$:

$$R_1 = \frac{0.2853}{0.3386} = 0.8425$$

$$P_1 = \frac{1}{1 + 0.8425^2} + 0.1084 = 0.6933$$

Similarly, the performance scores for all alternatives are calculated as $P_1 = 0.6933$, $P_2 = 0.7401$, $P_3 = 0.2431$, $P_4 = 0.6377$, $P_5 = 0.5244$. Accordingly, the ranking is A2 > A1 > A4 > A5 > A3, with A2 being the top-ranked alternative by PROBID.

The **advantages of PROBID** are as follows. (1) It uses several tiers of positive and negative ideal solutions, mitigating the impact of outliers, as no single outlier can dominate the references (Park et al., 2023). (2) It balances out extremes by factoring in the centroid (average) of all alternatives. (3) Rank reversal is mitigated due to the use of multiple reference points, which can provide a more comprehensive evaluation framework.

The **limitations of PROBID** are as follows. (1) Computing distances to multiple tiers of ideal solutions and an average solution is computationally intensive compared to other reference-type methods, making it potentially less efficient when handling a large ACM. (2) It is not completely outlier-proof, as extreme outliers can still affect the formation of different ideal tiers and/or distort the average solution.

## 7.11   Summary

Main points in this chapter of reference-type MCDM methods are as follows.

1. This chapter describes and illustrates 9 reference-type MCDM methods, each with distinct approaches to defining reference solutions and computing performance scores. Listed chronologically, these methods are TOPSIS, GRA, VIKOR, EDAS, MABAC, CODAS, PIV, MARCOS, and PROBID.
2. The methods covered are chosen based on either their widespread adoption (e.g., TOPSIS and VIKOR) or their recent development (e.g., MARCOS and PROBID).

3. Reference-type MCDM methods rank alternatives by comparing them to specific reference points (or solutions such as PIS and NIS) derived from the ACM. These reference points serve as benchmarks for evaluating the performance of alternatives.

4. Each method is clearly presented with its theoretical foundations, step-by-step algorithm, and numerical calculations to enhance reader comprehension. The calculations are performed using a common ACM dataset, allowing for easy comparison across different methods.

5. Each method has its own advantages and limitations, as presented in its respective section.

6. Rank reversal is a common limitation, as adding or removing alternatives can modify the reference solutions. As a result, the relative rankings of existing alternatives may shift based on the updated reference points.

7. Some methods primarily rely on distances to positive and/or negative ideal solutions (e.g., TOPSIS), while some are based on distances to an arithmetic or geometric mean solution (e.g., MABAC). Additionally, some methods incorporate both ideal and mean solutions (e.g., PROBID).

8. Understanding the principles of each method can help decision-makers choose the most appropriate one based on the specific decision-making context.

9. Mastering these reference-type MCDM methods allows decision-makers, researchers, and practitioners to effectively apply them to engineering and scientific applications that require systematic decision analysis.

10. Finally, as anticipated, the ranking results from the 9 MCDM methods applied to the illustrated example show variation in the top-ranked alternative. While some methods consistently prioritize the same alternative, others yield different rankings due to their distinct principles and ranking mechanisms. For the illustrated example, A2 is the top-ranked alternative according to TOPSIS, CODAS, MARCOS, and PROBID; A4 is ranked highest by GRA, EDAS, MABAC, and PIV; and A1 is the top-ranked alternative under VIKOR. Generally, while A2 and A4 appear frequently at the top, the ranking order of other alternatives varies across methods, as seen in each individual section. These ranking results cannot be generalized, and it is better to try several MCDM methods for any application.

## 7.13 Exercises

E7.1 You are given an ACM consisting of 4 alternatives (A1, A2, A3, A4) and 3 criteria (C1, C2, C3), with the following values:

| Alternatives | C1 | C2 | C3 |
|---|---|---|---|
| A1 | 0.93 | 600 | 8.25 |
| A2 | 0.51 | 700 | 6.33 |
| A3 | 0.77 | 500 | 3.16 |
| A4 | 0.82 | 400 | 2.98 |

Criteria C1 and C2 are the benefit criteria to be maximized, whereas C3 is a cost criterion to be minimized. The weights for the criteria are: $w_1 = 0.5$, $w_2 = 0.3$, and $w_3 = 0.2$. Apply all (or some of) the reference-type MCDM methods covered in this chapter to rank the alternatives. Identify the top-ranked alternative by each method and compare them.

E7.2 You are given an ACM consisting of 5 alternatives (A1, A2, A3, A4, A5) and 4 criteria (C1, C2, C3, C4), with the following values:

| Alternatives | C1 | C2 | C3 | C4 |
|---|---|---|---|---|
| A1 | 234 | 0.122 | 90.3 | 0.069 |
| A2 | 179 | 0.641 | 13.2 | 0.032 |
| A3 | 398 | 0.782 | 67.1 | 0.191 |
| A4 | 273 | 0.979 | 49.8 | 0.264 |
| A5 | 278 | 0.543 | 86.8 | 0.219 |

Criteria C2 and C3 are the benefit criteria to be maximized, whereas C1 and C4 are the cost criteria to be minimized. The weights for the criteria are: $w_1 = 0.3$, $w_2 = 0.4$, $w_3 = 0.2$, and $w_4 = 0.1$. Apply all (or some of) the reference-type MCDM methods covered in this chapter to rank the alternatives. Identify the top-ranked alternative by each method and compare them. *Note that calculations in this and subsequent examples increase with the increasing number of alternatives and criteria. Hence, readers are advised to implement the calculations of each method in a Microsoft Excel worksheet.*

E7.3 You are given an ACM consisting of 7 alternatives (A1, A2, A3, A4, A5, A6, A7) and 5 criteria (C1, C2, C3, C4, C5), with the following values:

| Alternatives | C1 | C2 | C3 | C4 | C5 |
|---|---|---|---|---|---|
| A1 | 3405 | 87.4 | 0.245 | 0.105 | 4.2 |
| A2 | 2159 | 45.2 | 0.521 | 0.187 | 3.7 |
| A3 | 4782 | 72.1 | 0.684 | 0.274 | 5.1 |
| A4 | 3594 | 33.9 | 0.319 | 0.143 | 2.9 |
| A5 | 2911 | 94.3 | 0.753 | 0.238 | 4.8 |
| A6 | 4100 | 59.7 | 0.602 | 0.194 | 3.4 |
| A7 | 3317 | 80.6 | 0.438 | 0.165 | 4.7 |

Criteria C1, C2, and C3 are the benefit criteria to be maximized, whereas C4 and C5 are the cost criteria to be minimized. The weights for the criteria are: $w_1 = 0.25$, $w_2 = 0.2$, $w_3 = 0.3$, $w_4 = 0.15$, and $w_5 = 0.1$. Apply all (or some of) the reference-type MCDM methods covered in this chapter to rank the alternatives. Identify the top-ranked alternative by each method and compare them.

E7.4 You are given an ACM consisting of 10 alternatives (A1, A2, A3, A4, A5, A6, A7, A8, A9, A10) and 6 criteria (C1, C2, C3, C4, C5, C6), with the following values:

| Alternatives | C1 | C2 | C3 | C4 | C5 | C6 |
|---|---|---|---|---|---|---|
| A1 | 575 | 0.125 | 63.8 | 0.0215 | 12.3 | 3.5 |
| A2 | 432 | 0.315 | 89.2 | 0.0382 | 15.7 | 2.9 |
| A3 | 689 | 0.498 | 74.5 | 0.0497 | 18.2 | 4.1 |
| A4 | 540 | 0.276 | 95.3 | 0.0328 | 14.9 | 3.8 |
| A5 | 478 | 0.605 | 70.4 | 0.0409 | 16.5 | 4.5 |
| A6 | 615 | 0.451 | 82.6 | 0.0462 | 17.9 | 3.7 |
| A7 | 503 | 0.333 | 88.9 | 0.0375 | 13.8 | 3.2 |
| A8 | 389 | 0.254 | 59.1 | 0.0293 | 10.5 | 2.8 |
| A9 | 455 | 0.394 | 76.3 | 0.0319 | 14.2 | 4.3 |
| A10 | 612 | 0.512 | 81.9 | 0.0415 | 15.1 | 3.6 |

Criteria C1 and C4 are the benefit criteria to be maximized, whereas C2, C3, C5, C6 are the cost criteria to be minimized. The weights for the criteria are: $w_1 = 0.2$, $w_2 = 0.15$, $w_3 = 0.25$, $w_4 = 0.15$, $w_5 = 0.15$, and $w_6 = 0.1$. Apply all (or some of) the reference-type MCDM methods covered in this chapter to rank the alternatives. Identify the top-ranked alternative by each method and compare them.

E7.5 Describe the similarities and differences between TOPSIS and PIV.

E7.6 Describe the similarities and differences between EDAS and CODAS.

E7.7 Describe the similarities and differences between MABAC and PROBID.